\documentclass[a4paper,12pt]{article}

\usepackage{amsmath}
\usepackage{amsthm}
\usepackage{amssymb}
\usepackage{latexsym}
\usepackage{palatino}
\usepackage{pstricks}
\usepackage{graphicx, pst-plot, pst-node, pst-text, pst-tree}
\usepackage{cite}

\theoremstyle{plain}
\newtheorem{theorem}{Theorem}[section]
\newtheorem{corollary}[theorem]{Corollary}
\newtheorem{lemma}[theorem]{Lemma}
\newtheorem{proposition}[theorem]{Proposition}

\newtheorem{question}[theorem]{Question}
\theoremstyle{definition}

\theoremstyle{remark}

\newenvironment{proofsketch}{\noindent{\it Sketch of proof.}}{\qed\bigskip}

\newcounter{todocounter}

\makeatletter
\newfont{\footsc}{cmcsc10 at 8truept}
\newfont{\footbf}{cmbx10 at 8truept}
\newfont{\footrm}{cmr10 at 10truept}
\title{A Survey of Simple Permutations}
\author{Robert Brignall\\[-5pt]
\small Department of Mathematics\\[-5pt]
\small University of Bristol\\[-5pt]
\small Bristol, UK\\[-5pt]
\small \texttt{Robert.Brignall@bristol.ac.uk}\\[-5pt]
\small \texttt{http://www.maths.bris.ac.uk/\~{}marlfb}}

\date{\today \\[6pt]
}

\begin{document}
\maketitle

\newcommand{\Sub}{\operatorname{Sub}}
\newcommand{\Av}{\operatorname{Av}}
\newcommand{\rect}{\mathrm{rect}}
\newcommand{\A}{\mathcal{A}}
\newcommand{\B}{\mathcal{B}}
\newcommand{\C}{\mathcal{C}}
\newcommand{\D}{\mathcal{D}}
\renewcommand{\P}{\mathcal{P}}
\newcommand{\Q}{\mathcal{Q}}
\newcommand{\R}{\mathcal{R}}
\renewcommand{\S}{\mathcal{S}}
\newcommand{\wrc}[1]{\langle #1 \rangle}
\newcommand{\Si}{\operatorname{Si}}
\newcommand{\dom}{\operatorname{dom}}

\newcommand{\OEIS}[1]{(sequence #1 of~\cite{sloane:the-on-line-enc:})}

\begin{abstract}
We survey the known results about simple permutations. In particular, we present a number of recent enumerative and structural results pertaining to simple permutations, and show how simple permutations play an important role in the study of permutation classes. We demonstrate how classes containing only finitely many simple permutations satisfy a number of special properties relating to enumeration, partial well-order and the property of being finitely based.
\end{abstract}

\section{Introduction}

An \emph{interval} of a permutation $\pi$ corresponds to a set of contiguous indices $I=[a,b]$ such that the set of values $\pi(I)=\{\pi(i) : i\in I\}$ is also contiguous. Every permutation of length $n$ has intervals of lengths $0$, $1$ and $n$. If a permutation $\pi$ has no other intervals, then $\pi$ is said to be \emph{simple}. For example, the permutation $\pi=28146357$ is not simple as witnessed by the non-trivial interval $4635$ ($=\pi(4)\pi(5)\pi(6)\pi(7)$), while $\sigma=51742683$ is simple.\footnote{Simplicity may be defined for all relational structures; such structures have variously been called \emph{prime} or \emph{indecomposable}.}

While intervals of permutations have applications in biomathematics, particularly to genetic algorithms and the matching of gene sequences (see Corteel, Louchard, and Pemantle~\cite{corteel:common-interval:} for extensive references), simple permutations form the ``building blocks'' of permutation classes and have thus received intensive study in recent years. We will see in Section~\ref{sec-classes} the various ways in which simplicity plays a role in the study of permutation classes, but we begin this short survey by introducing the substitution decomposition in Subsection~\ref{subsec-subst-decomp}, and thence by reviewing the structural and enumerative results of simple permutations themselves in Section~\ref{sec-structure}. The rest of this subsection will cover several basic definitions that we will require.

Two finite sequences of the same length, $\alpha=a_1a_2\cdots a_n$ and $\beta=b_1b_2\cdots b_n$, are said to be \emph{order isomorphic} if, for all $i,j$, we have $a_i < a_j$ if and only if $b_i < b_j$. As such, each sequence of distinct real numbers is order isomorphic to a unique permutation. Similarly, any given subsequence (or \emph{pattern}) of a permutation $\pi$ is order isomorphic to a smaller permutation, $\sigma$ say, and such a subsequence is called a \emph{copy} of $\sigma$ in $\pi$.  We may also say that $\pi$ \emph{contains} $\sigma$ (or, in some texts, $\pi$ \emph{involves} $\sigma$) and write $\sigma \leq \pi$. If, on the other hand, $\pi$ does not contain a copy of some given $\sigma$, then $\pi$ is said to \emph{avoid} $\sigma$.  For example, $\pi=918572346$ contains $51342$ because of the subsequence $91572$ ($=\pi(1)\pi(2)\pi(4)\pi(5)\pi(6)$), but avoids $3142$.

\begin{figure}
\begin{center}
\psset{xunit=0.01in, yunit=0.01in}
\psset{linewidth=0.005in}
\begin{pspicture}(0,0)(90,90)
\psaxes[dy=10,Dy=1,dx=10,Dx=1,tickstyle=bottom,showorigin=false,labels=none](0,0)(90,90)
\pscircle*(10,90){0.04in}
\pscircle*(20,30){0.04in}
\pscircle*(30,40){0.04in}
\pscircle*(40,80){0.04in}
\pscircle*(50,20){0.04in}
\pscircle*(60,60){0.04in}
\pscircle*(70,70){0.04in}
\pscircle*(80,10){0.04in}
\pscircle*(90,50){0.04in}
\end{pspicture}
\end{center}
\caption{The plot of the permutation $\pi=934826715$.}\label{fig-perm-ex}
\end{figure}
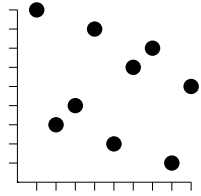

It will often be useful to view permutations and order isomorphism graphically. Two sets $S$ and $T$ of points in the plane are said to be order isomorphic if the axes for the set $S$ can be stretched and shrunk in some manner to map the points of $S$ bijectively onto the points of $T$, i.e., if there are strictly increasing functions $f,g:\mathbb{R}\rightarrow\mathbb{R}$ such that $\{(f(s_1),g(s_2)) : (s_1,s_2)\in S\}=T$. Note that this forms an equivalence relation since the inverse of a strictly increasing function is also strictly increasing. The {\it plot\/} of the permutation $\pi$ is then the point set $\{(i,\pi(i))\}$, and every finite point set in the plane in which no two points share a coordinate (often called a {\it generic\/} or {\it noncorectilinear\/} set) is order isomorphic to the plot of a unique permutation (see Figure~\ref{fig-perm-ex} for an example). Note that, with a slight abuse of notation, we will say that a point set is order isomorphic to a permutation.

The pattern containment order forms a partial order on the set of all permutations. Downsets of permutations under this order are called \emph{permutation classes}.\footnote{In the past, permutation classes have also been called \emph{closed classes} or \emph{pattern classes}.} In other words, if $\C$ is a permutation class and $\pi\in \C$, then for any permutation $\sigma$ with $\sigma\leq\pi$ we have $\sigma\in\C$. A given permutation class is often described in terms of its minimal avoidance set, or \emph{basis}. More formally, the basis $B$ of a permutation class $\C$ is the smallest set for which $\C = \{\pi\mid\beta\not\leq\pi\textrm{ for all }\beta\in B\}$. For a permutation class $\C$, we denote by $\C_n$ the set $\C \cap S_n$, i.e.\ the permutations in $\C$ of length $n$, and we refer to $f(x) = \sum |\C_n| x^n$ as the {\it generating function for $\C$\/}.

Analogues of pattern containment exist for other relational structures; sets of structures closed under taking induced substructures are known as \emph{hereditary properties}. Hereditary properties of graphs have received considerable attention (see Bollob{\'a}s~\cite{bollobas:hereditary-prop:} for a survey of some older results), while more recently attention has been given to hereditary properties of a variety of structures including tournaments, ordered graphs and posets (see, for example, Balogh {\it et al.\/}~\cite{balogh:hereditary-prop:,balogh:hereditary-prop:a,balogh:hereditary-prop:b}, and Bollob{\'a}s's recent survey~\cite{bollobas:hereditary-and:}).

\subsection{Substitution Decomposition}\label{subsec-subst-decomp}

The simple permutations form the elemental building blocks upon which all other permutations are constructed by means of the substitution decomposition.\footnote{This decomposition is also called the modular decomposition, disjunctive decomposition or $X$-join in other contexts.} Analogues of this decomposition exist for every relational structure, and it has frequently arisen in a wide variety of perspectives, ranging from game theory to combinatorial optimization --- for references see M\"ohring~\cite{mohring:algorithmic-asp:a} or M\"ohring and Radermacher~\cite{mohring:substitution-de:}. Its first appearance seems to be in a 1953 talk by Fra{\"{\i}}ss{\'e} (though only the abstract of this talk~\cite{fraisse:on-a-decomposit:} survives). It did not appear in an article until Gallai~\cite{gallai:transitiv-orien:} (for an English translation, see~\cite{gallai:a-translation-o:}), who applied them particularly to the study of transitive orientations of graphs.

Given a permutation $\sigma$ of length $m$ and nonempty permutations $\alpha_1,\dots,\alpha_m$, the {\it inflation\/} of $\sigma$ by $\alpha_1,\dots,\alpha_m$ --- denoted $\sigma[\alpha_1,\dots,\alpha_m]$ --- is the permutation obtained by replacing each entry $\sigma(i)$ by an interval that is order isomorphic to $\alpha_i$.  For example, $2413[1,132,321,12]=479832156$.  Conversely, a
\emph{deflation} of $\pi$ is any expression of $\pi$ as an inflation $\pi = \sigma[\pi_1,\pi_2,\ldots,\pi_m]$, and we will call $\sigma$ a \emph{skeleton} of $\pi$. We then have the \emph{substitution decomposition} of permutations:

\begin{proposition}[Albert and Atkinson~\cite{albert:simple-permutat:}]\label{simple-decomp-2}
Every permutation may be written as the inflation of a unique simple permutation. Moreover, if $\pi$ can be written as $\sigma[\alpha_1,\dots,\alpha_m]$ where $\sigma$ is simple and $m\ge 4$, then the $\alpha_i$'s are unique.
\end{proposition}

Non-unique cases arise when a permutation can be written as an inflation of either $12$ or $21$, and to recover uniqueness we may choose a particular decomposition in a variety of ways. The one we will use is as follows.

\begin{proposition}[Albert and Atkinson~\cite{albert:simple-permutat:}]\label{simple-decomp-3}
If $\pi$ is an inflation of $12$, then there is a unique sum indecomposable $\alpha_1$ such that $\pi=12[\alpha_1,\alpha_2]$ for some $\alpha_2$, which is itself unique.  The same holds with $12$ replaced by $21$ and ``sum'' replaced by ``skew''.
\end{proposition}

\begin{figure}
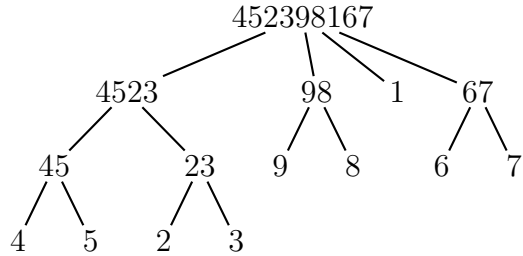

\begin{center}
\pstree[nodesep=2pt,levelsep=5ex]{\TR{$452398167$} }{
  \pstree{ \TR{$4523$} }{
    \pstree{ \TR{$45$} }{
      \TR{$4$}
      \TR{$5$}
    }
    \pstree{ \TR{$23$} }{
      \TR{$2$}
      \TR{$3$}
    }
  }
  \pstree{ \TR{$98$} }{
    \TR{$9$}
    \TR{$8$}
  }
  \TR{$1$}
  \pstree{ \TR{$67$} }{
    \TR{$6$}
    \TR{$7$}
  }
}
\end{center}
\caption{The substitution decomposition tree of $\pi=452398167$.}\label{fig-substdecomptree-1}
\end{figure}

The substitution decomposition tree for a permutation is obtained by recursively decomposing until we are left only with inflations of simple permutations by singletons. For example, the permutation $\pi=452398167$ is decomposed as
\begin{eqnarray*}
452398167 &=& 2413[3412,21,1,12]\\
 &=& 2413[21[12,12],21[1,1],1,12[1,1]]\\
 &=& 2413[21[12[1,1],12[1,1]],21[1,1],1,12[1,1]]
\end{eqnarray*}
and its substitution decomposition tree is given in Figure~\ref{fig-substdecomptree-1}.

\paragraph{Computation in Linear Time.}
The substitution decomposition is most frequently used in solving algorithmic problems, and consequently much attention has been given to its computation in optimal time.\footnote{In particular, graph decomposition has received significant attention, with the first $O(|V| + |E|)$ algorithms appearing in 1994 by McConnell and Spinrad~\cite{mcconnell:linear-time-mod:} and Cournier and Habib~\cite{cournier:a-new-linear-al:}.} By its connection to the intervals of a permutation, a first approach to compute the substitution decomposition might be simply to compute all the intervals of our given permutation. Since there may be as many as $N=n(n-1)/2$ such intervals in a permutation of length $n$, listing these will not yield a linear $O(n)$ algorithm for the substitution decomposition. However, this computation has received significant attention through its connections with biomathematics, with an $O(n+N)$ time algorithm being given by Bergeron, Chauve, Montgolfier and Raffinot~\cite{bergeron:computing-commo:}.\footnote{In fact, Bergeron {\it\/et al\/} show how to compute the ``common intervals'' --- a generalisation of our notion of interval applied to sets of permutations.}

The first algorithm to compute the substitution decomposition of a permutation in linear time was given by Uno and Yagiura~\cite{uno:fast-algorithms:}, while Bergeron {\it et al} have since given a simpler algorithm. A \emph{strong interval} of a permutation $\pi$ is an interval $I$ for which every other interval $J$ satisfies one of $J\subseteq I$, $I\subseteq J$ or $I\cap J=\emptyset$. For example, given $\pi=234615$, the interval $234$ ($=\pi(1)\pi(2)\pi(3)$) is a strong interval, but $23$ is not, because it has non-trivial intersection with $34$. A permutation can have at most $2n-1$ strong intervals (note that the $n$ singletons and the whole permutation are all strong intervals), and Bergeron {\it et al\/} give an optimal $O(n)$ algorithm to list them all. The substitution decomposition tree of the permutation follows immediately.

It is worth noticing that this algorithm does not give the simple skeletons for each internal node of the decomposition tree --- indeed, there are currently no linear time algorithms to do this. It is, however, straightforward to compute the label of any particular node in linear time, e.g. by finding a representative symbol for each strong interval lying below the node, and then computing the permutation order isomorphic to this sequence of representatives.

%
%
%
%
%
%
%
%
%
%
\section{Enumeration and Structure}\label{sec-structure}
\subsection{Enumeration and Asymptotics}
The number of simple permutations of length $n=1,2,\ldots$ is $s_n=1,2,0,2,6,46,338,2926,28146,\ldots$ \OEIS{A111111}. Albert, Atkinson and Klazar~\cite{albert:the-enumeration:} showed, in a straightforward argument making use of the substitution decomposition, that the sequence $(s_n)_{n\geq 4}$ is given by \[s_n = - \mathrm{Com}_n + (-1)^{n+1}\cdot 2,\]
where $\mathrm{Com}_n$ is the coefficient of $x^n$ in the functional inverse of $f(x)=\sum_{n=1}^\infty n!x^n$ \OEIS{A059372}.\footnote{The term $\mathrm{Com}_n$ is used since the function $f^{-1}(x)$ first appeared in an exercise in Comtet~\cite{comtet:advanced-combinat:}.}

Asymptotically, the sequence $s_n$ may be counted using a probabilistic argument based on counting the intervals of a random permutation. Let the random variable $X_k$ denote the number of intervals of length $k$ in a random permutation $\pi$ of length $n$. An interval of length $k$ may be viewed as a mapping from a contiguous set of positions to a contiguous set of values, for which the set of positions must begin at one of the first $n-k+1$ positions of $\pi$, and the lowest point in the set of values must be one of the lowest $n-k+1$ values of $\pi$. Of the $\binom{n}{k}$ sets of values to which the contiguous set of positions may be mapped, only one maps to the chosen contiguous set of values. Thus we have
\[ \mathbb{E}[X_k] = \frac{(n-k+1)^2}{\binom{n}{k}} = \frac{(n-k+1)(n-k+1)!k!}{n!}.\]
Our first observation is that, as $n\rightarrow \infty$, $\mathbb{E}[X_2] = \frac{2(n-1)}{n}\rightarrow 2$. Thus, asymptotically, we should expect to find precisely two intervals of size two in a random permutation. We are seeking the asymptotics of the expected number of proper intervals, i.e.\ the sum $\displaystyle\sum_{k=2}^{n-1}\mathbb{E}[X_k]$, and want to demonstrate that $\displaystyle\sum_{k=3}^{n-1}\mathbb{E}[X_k]\rightarrow 0$ as $n\rightarrow 0$. We first consider the cases $k=3$, $k=4$, $k=n-2$ (assuming $n\geq 4$) and $k=n-1$ separately:
\begin{eqnarray*}
\mathbb{E}[X_3] &=& \frac{6(n-2)}{n(n-1)} \leq \frac{6}{n} \rightarrow 0\\
\mathbb{E}[X_4] &=& \frac{4!(n-3)}{n(n-1)(n-2)} \leq \frac{24}{n^2} \rightarrow 0\\
\mathbb{E}[X_{n-2}] &=& \frac{3\cdot 3!}{n(n-1)} \leq \frac{24}{n^2} \rightarrow 0\\
\mathbb{E}[X_{n-1}] &=& \frac{4}{n} \rightarrow 0.
\end{eqnarray*}
The remaining terms form a partial sum, which converges providing $\displaystyle \frac{\mathbb{E}[X_{k+1}]}{\mathbb{E}[X_k]}<1$. Simplifying this equation gives $2k^2-(3n+1)k + n^2+n+1 >0$, a quadratic in $k$, which yields two roots. The smaller of these satisfies $0< k^-\leq n$, the larger $k^+>n$. Thus for $k\leq k^-$, $\mathbb{E}[X_k]$ is decreasing, while for $k^-<k<n$, $\mathbb{E}[X_k]$ is increasing, and hence $\mathbb{E}[X_k]\leq 24/n^2$ for $4\leq k\leq n-2$. Thus
\[\sum_{k=4}^{n-2} \mathbb{E}[X_k] \leq (n-5)\frac{24}{n^2} \leq \frac{24}{n} \rightarrow 0.\]
Subsequently, the only term of $\displaystyle\sum_{k=2}^{n-1}\mathbb{E}[X_k]$ which is non-zero in the limit $n\rightarrow\infty$ is $k=2$.\footnote{A similar argument can be applied to graphs, but in this case we find that $\mathbb{E}[X_k]\rightarrow 0$ as $n\rightarrow 0$ for every $2\leq k\leq n-1$, and thus, asymptotically, almost all graphs are indecomposable. The same applies to tournaments, posets, and (more generally) structures defined on a single asymmetric relation --- see M\"ohring \cite{mohring:on-the-distribu:}.}

Ignoring larger intervals, occurrences of intervals of size $2$ in a large random permutation $\pi$ can roughly be regarded as independent events, and as we know the expectation of $X_2$ is $2$, the occurrence of any specific interval is relatively rare. Heuristically, this suggests that $X_2$ is asymptotically Poisson distributed with parameter $2$. Using this heuristic, we have $\textrm{Pr}(X_2=0)\rightarrow e^{-2}$ as $n\rightarrow\infty$, and so there are approximately $\frac{n!}{e^2}$ simple permutations of length $n$.

A formal argument for this was implicitly given by Uno and Yagiura~\cite{uno:fast-algorithms:}, and was made explicit by Corteel, Louchard, and Pemantle~\cite{corteel:common-interval:}. The method, however, essentially dates back to the 1940s with Kaplansky~\cite{kaplansky:the-asymptotic-:} and Wolfowitz~\cite{wolfowitz:note-on-runs-of:}, who considered ``runs'' within permutations --- a \emph{run} is a set of points with contiguous positions whose values are $i, i+1,\ldots, i+r$ or $i+r, i+r-1,\ldots,i$, in that order.\footnote{Atkinson and Stitt~\cite{atkinson:restricted-perm:a} called permutations containing no runs \emph{strongly irreducible}. Note that this is equivalent to a permutation containing no intervals of size two.} A non-probabilistic approach for these first order asymptotics based on Lagrange inversion can be obtained from a more general theorem of Bender and Richmond~\cite{bender:an-asymptotic-expansion:}.

More precise asymptotics, meanwhile, have been found using a non-probabilistic method (but one relying on the work of Kaplansky) by Albert, Atkinson, and Klazar~\cite{albert:the-enumeration:}. They obtain the following theorem, and note that higher order terms are calculable given sufficient computation:

\begin{theorem}[Albert, Atkinson and Klazar~\cite{albert:the-enumeration:}]
The number of simple permutations of length $n$ is asymptotically given by \[\frac{n!}{e^2}\left(1-\frac{4}{n}+\frac{2}{n(n-1)}+O(n^{-3})\right).\]
\end{theorem}

\subsection{Exceptional Simple Permutations}
Given a simple permutation $\pi$, one might ask what simple permutations are contained within $\pi$. In particular, is there a point that can be removed from $\pi$ to leave a sequence order isomorphic to a simple permutation? This is not quite true, but allowing one- or two-point deletions suffices. The following theorem is a special case of a more general result for every relational structure whose relations are binary and irreflexive.

\begin{theorem}[Schmerl and Trotter~\cite{schmerl:critically-inde:}]\label{thm-schmerl-trotter}
Every simple permutation of length $n\ge 2$ contains a simple permutation of length $n-1$ or $n-2$.
\end{theorem}

\begin{figure}
\begin{center}
\begin{tabular}{ccccccc}
\psset{xunit=0.01in, yunit=0.01in}
\psset{linewidth=0.005in}
\begin{pspicture}(0,0)(100,100)
\psaxes[dy=10, Dy=1, dx=10, Dx=1, tickstyle=bottom, showorigin=false, labels=none](0,0)(100,100)
\pscircle*(10,60){0.04in}
\pscircle*(20,50){0.04in}
\pscircle*(30,70){0.04in}
\pscircle*(40,40){0.04in}
\pscircle*(50,80){0.04in}
\pscircle*(60,30){0.04in}
\pscircle*(70,90){0.04in}
\pscircle*(80,20){0.04in}
\pscircle*(90,100){0.04in}
\pscircle*(100,10){0.04in}
\end{pspicture}
&\rule{10pt}{0pt}&
\psset{xunit=0.01in, yunit=0.01in}
\psset{linewidth=0.005in}
\begin{pspicture}(0,0)(100,100)
\psaxes[dy=10,Dy=1,dx=10,Dx=1,tickstyle=bottom,showorigin=false,labels=none](0,0)(100,100)
\pscircle*(10,100){0.04in}
\pscircle*(20,80){0.04in}
\pscircle*(30,60){0.04in}
\pscircle*(40,40){0.04in}
\pscircle*(50,20){0.04in}
\pscircle*(60,10){0.04in}
\pscircle*(70,30){0.04in}
\pscircle*(80,50){0.04in}
\pscircle*(90,70){0.04in}
\pscircle*(100,90){0.04in}
\end{pspicture}
&\rule{10pt}{0pt}&
\psset{xunit=0.01in, yunit=0.01in}
\psset{linewidth=0.005in}
\begin{pspicture}(0,0)(100,100)
\psaxes[dy=10,Dy=1,dx=10,Dx=1,tickstyle=bottom,showorigin=false,labels=none](0,0)(100,100)
\pscircle*(10,10){0.04in}
\pscircle*(20,30){0.04in}
\pscircle*(30,50){0.04in}
\pscircle*(40,70){0.04in}
\pscircle*(50,90){0.04in}
\pscircle*(60,20){0.04in}
\pscircle*(70,40){0.04in}
\pscircle*(80,60){0.04in}
\pscircle*(90,80){0.04in}
\pscircle*(100,100){0.04in}
\end{pspicture}
&\rule{10pt}{0pt}&
\psset{xunit=0.01in, yunit=0.01in}
\psset{linewidth=0.005in}
\begin{pspicture}(0,0)(100,100)
\psaxes[dy=10, Dy=1, dx=10, Dx=1, tickstyle=bottom, showorigin=false, labels=none](0,0)(100,100)
\pscircle*(10,50){0.04in}
\pscircle*(20,100){0.04in}
\pscircle*(30,40){0.04in}
\pscircle*(40,90){0.04in}
\pscircle*(50,30){0.04in}
\pscircle*(60,80){0.04in}
\pscircle*(70,20){0.04in}
\pscircle*(80,70){0.04in}
\pscircle*(90,10){0.04in}
\pscircle*(100,60){0.04in}
\end{pspicture}
\end{tabular}
\end{center}
\caption{The two permutations on the left are wedge alternations, the two on the right are parallel alternations.}
\label{fig-wedge-parallel}
\end{figure}
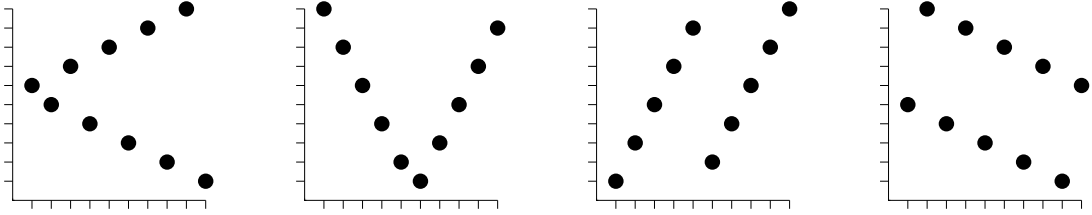

In most cases, however, a single point deletion is sufficient. If none of the one point deletions of a given simple permutation $\pi$ is simple, then $\pi$ is said to be \emph{exceptional}. Schmerl and Trotter call such structures \emph{critically indecomposable}, and present a complete characterisation in the analogous problem for partially ordered sets.

\begin{figure}
\begin{center}
\begin{tabular}{ccc}
\psset{xunit=0.01in, yunit=0.01in}
\psset{linewidth=0.005in}
\begin{pspicture}(0,0)(120,120)
\psaxes[dy=10, Dy=1, dx=10, Dx=1, tickstyle=bottom, showorigin=false,labels=none](0,0)(120,120)
\pscircle*(10,70){0.04in}
\pscircle*(20,50){0.04in}
\pscircle*(30,80){0.04in}
\pscircle*(40,40){0.04in}
\pscircle*(50,90){0.04in}
\pscircle*(60,30){0.04in}
\pscircle*(70,100){0.04in}
\pscircle*(80,20){0.04in}
\pscircle*(90,110){0.04in}
\pscircle*(100,10){0.04in}
\pscircle*(110,120){0.04in}
\pscircle*(120,60){0.04in}
\end{pspicture}
&\rule{10pt}{0pt}&
\psset{xunit=0.01in, yunit=0.01in}
\psset{linewidth=0.005in}
\begin{pspicture}(0,0)(120,120)
\psaxes[dy=10, Dy=1, dx=10, Dx=1, tickstyle=bottom, showorigin=false, labels=none](0,0)(120,120)
\pscircle*(10,60){0.04in}
\pscircle*(20,120){0.04in}
\pscircle*(30,50){0.04in}
\pscircle*(40,70){0.04in}
\pscircle*(50,40){0.04in}
\pscircle*(60,80){0.04in}
\pscircle*(70,30){0.04in}
\pscircle*(80,90){0.04in}
\pscircle*(90,20){0.04in}
\pscircle*(100,100){0.04in}
\pscircle*(110,10){0.04in}
\pscircle*(120,110){0.04in}
\end{pspicture}
\end{tabular}
\end{center}
\caption{The two types of wedge simple permutation, type $1$ (left) and type $2$ (right).}\label{fig-fund-wedge}
\end{figure}
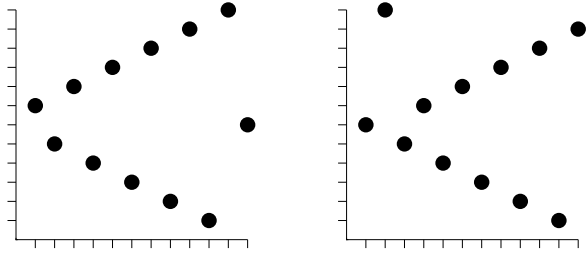

To consider the exceptional simple permutations, we first define a set of permutations called \emph{alternations}. A {\it horizontal alternation\/} is a permutation in which every odd entry lies to the left of every even entry, or the reverse of such a permutation.  Similarly, a {\it vertical alternation\/} is the group-theoretic inverse of a horizontal alternation. Of these alternations, we identify two families in which each ``side'' of the alternation forms a monotone sequence, namely the {\it parallel\/} and {\it wedge alternations\/} --- see Figure~\ref{fig-wedge-parallel} for definitions. While any parallel alternation is already simple or very nearly so, wedge alternations are not. Any wedge alternation may be extended to form a simple permutation by placing a single point in one of two places, thus forming wedge simple permutations of types 1 and 2 --- see Figure~\ref{fig-fund-wedge}.

The exceptional simple permutations turn out to be precisely the set of parallel alternations:

\begin{theorem}[Albert and Atkinson~\cite{albert:simple-permutat:}]
The only simple permutations that do not have a one point deletion are the simple parallel alternations, i.e.\ those of the form
\[ 246\cdots (2m)135\cdots (2m-1) \quad (m\ge 2)\]and every symmetry of this permutation.
\end{theorem}

\begin{proof}
Define a poset $(P_\pi,\prec_\pi)$ of the permutation $\pi$ by $x\prec_\pi y$ if and only if $x<y$ and $\pi(x)<\pi(y)$. Note that the poset $(P_\pi,\prec_\pi)$ has dimension $2$, and conversely that all posets of dimension $2$ correspond to a unique permutation, up to permutation inverses. A permutation is simple if and only if its corresponding poset is also simple, i.e.\ $P_\pi$ contains no proper nonsingleton subset $I$ for which every pair $x$, $y$ of $I$ are ordered with respect to all elements of $P_\pi\setminus I$ in the same way. Furthermore, $\pi$ is exceptional if and only if the corresponding poset is also exceptional (or, more usually, critically indecomposable). Schmerl and Trotter~\cite{schmerl:critically-inde:} classify all the critically indecomposable posets, from which the result follows.
\end{proof}

\subsection{Pin Sequences and Decomposition}

In the plot of a permutation, an interval is clearly identified as a set of points enclosed in an axes-parallel rectangle, with no points lying in the regions directly above, below, to the left or to the right. Conversely, since simple permutations contain no non-trivial intervals, any axes-parallel rectangle drawn over the plot of a simple permutation must be ``separated'' by at least one point in one of these four regions (unless the rectangle contains only one point or the whole permutation). Extending our rectangle to include this extra point puts us in a similar situation, and inductively our rectangle will eventually be extended to contain every point of the permutation. This is the motivating idea behind ``pin sequences''.

Given points $p_1,\ldots,p_m$ in the plane, we denote by $\rect(p_1,\ldots,p_m)$ the smallest axes-parallel rectangle containing them. A \emph{pin} for the points $p_1,\ldots,p_m$ is any point $p_{m+1}$ not contained in $\rect(p_1,\ldots,p_m)$ that lies either horizontally or vertically amongst them. Such a pin will have a \emph{direction} --- up, down, left or right --- which records the position of $p_{m+1}$ relative to $\rect(p_1,\ldots,p_m)$. A \emph{pin sequence} is a sequence of points $p_1,p_2,\ldots$ such that $p_i$ is a pin for the points $p_1,\ldots,p_{i-1}$ for every $i\geq 3$, while a \emph{proper pin sequence} is a pin sequence satisfying two further conditions:

\begin{itemize}
\item \emph{Maximality condition}: each pin must be maximal in its direction. That is, in the plot of a permutation, of all possible pins for $p_1,\ldots,p_{i-1}$ having the same direction as $p_i$ ($i\geq 3$), the pin $p_i$ is furthest from $\rect(p_1,\ldots,p_{i-1})$.\footnote{In certain situations, notably where pin sequences are used to construct permutations from scratch (rather than taking pin sequences from the points in the plot of a permutation), the maximality condition is replaced with the \emph{externality condition}, requiring that $p_{i+1}$ lies outside $\rect(p_1,\ldots,p_i)$.}
\item \emph{Separation condition}: $p_{i+1}$ must separate $p_i$ from $\rect(p_1,\ldots,p_{i-1})$ ($i\geq 2$). That is, $p_{i+1}$ must lie horizontally or vertically between $p_i$ and $\rect(p_1,\ldots,p_{i-1})$.
\end{itemize}

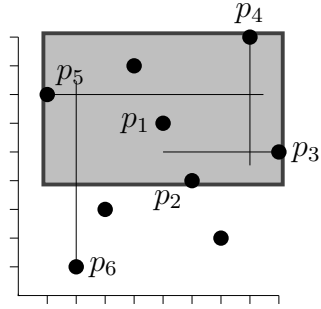
\begin{figure}
\begin{center}
\psset{xunit=0.01in, yunit=0.01in}
\psset{linewidth=0.005in}
\begin{pspicture}(0,0)(145,145)
\psaxes[dy=15,dx=15,tickstyle=bottom,showorigin=false,labels=none](0,0)(135,135)
\psframe[linecolor=darkgray,fillstyle=solid,fillcolor=lightgray,linewidth=0.02in](12,57)(138,138)
\psline(75,75)(135,75)
\psline(120,135)(120,68)
\psline(15,105)(127,105)
\psline(30,15)(30,112)
\pscircle*(15,105){0.04in}
\pscircle*(30,15){0.04in}
\pscircle*(45,45){0.04in}
\pscircle*(60,120){0.04in}
\pscircle*(75,90){0.04in}
\pscircle*(90,60){0.04in}
\pscircle*(105,30){0.04in}
\pscircle*(120,135){0.04in}
\pscircle*(135,75){0.04in}
\uput[l](75,90){$p_1$}
\uput[dl](90,60){$p_2$}
\uput[r](135,75){$p_3$}
\uput[u](120,135){$p_4$}
\uput[ur](15,105){$p_5$}
\uput[r](30,15){$p_6$}
\end{pspicture}
\end{center}
\caption{A pin sequence of the simple permutation $713864295$. The shaded box represents $\rect(p_1,p_2,p_3,p_4,p_5)$.}
\label{fig-pins-first}
\end{figure}

See Figure~\ref{fig-pins-first} for an example. Proper pin sequences are intimately connected with simple permutations. In one direction, we have:
\begin{theorem}[Brignall, Huczynska and Vatter~\cite{brignall:simple-permutat:a}]\label{pins-simple}
If $p_1,\dots,p_m$ is a proper pin sequence of length $m\geq 5$ then one of the sets of points $\{p_1,\dots,p_m\}$, $\{p_1,\dots,p_m\}\setminus\{p_1\}$, or $\{p_1,\dots,p_m\}\setminus\{p_2\}$ is order isomorphic to a simple permutation.
\end{theorem}

This should come as no surprise --- by our motivation, pin sequences encapsulate precisely what it means to be simple. Moreover, we also expect to be able to find long proper pin sequences within arbitrary simple permutations. If no such pin sequence exists, then we encounter the two families of alternations defined in the previous subsection.

\begin{theorem}[Brignall, Huczynska and Vatter~\cite{brignall:simple-permutat:a}]\label{sp2-really-main}
Every simple permutation of length at least $2(256k^8)^{2k}$ contains a proper pin sequence of length $2k$, a parallel alternation of length $2k$, or a wedge simple permutation of length $2k$.
\end{theorem}

\begin{proofsketch}
Suppose that a simple permutation $\pi$ of length $n$ contains neither a proper pin sequence of length $2k$ nor a parallel or wedge alternation of length $2k$. A pin sequence $p_1,p_2,\ldots,p_n$ of $\pi$ is said to be \emph{right-reaching} if $p_n$ corresponds to the rightmost point of $\pi$. Moreover, from each pair of points in a simple permutation $\pi$ there exists a right-reaching proper pin sequence. We now consider the collection of $\lfloor n/2\rfloor$ proper right-reaching pin sequences of $\pi$ beginning with the first and second points, the third and fourth points, and so on, reading from left to right.

Two pin sequences $p_1,p_2,\ldots$ and $q_1,q_2,\ldots$ are said to \emph{converge at the point $x$} if there exists $i$ and $j$ such that $p_i=q_j=x$ but $\{p_1,\ldots,p_{i-1}\}$ and $\{q_1,\ldots,q_{i-1}\}$ are disjoint. Our $\lfloor n/2\rfloor$ proper right-reaching pin sequences all end at the rightmost point of $\pi$, and so every pair of these must converge at some point.

If $16k$ proper pin sequences converge at a point $x$, then we may show there must exist an alternation of length at least $2k$. Such an alternation may not necessarily be parallel or wedge, but, by the Erd\H os-Szekeres Theorem, every alternation of length at least $2k^4$ must contain a parallel or wedge alternation of length at least $2k$. Accordingly, since $\pi$ contains no parallel or wedge alternation of length $2k$, fewer than $16k^4$ of our pin sequences can converge at any point. Since no pin sequence contains $2k$ or more pins, the number of right-reaching pin sequences is bounded by how many can converge at each steps towards the rightmost pin, and there can be at most $2k-1$ steps. Thus $\lfloor n/2\rfloor < 2(16k^4)^{2k}$.

Finally, if this process has produced a wedge alternation, it is then necessary to demonstrate the existence of one additional point to form a wedge simple permutation of type 1 or 2. For a wedge alternation containing $4k^2$ points and oriented $<$ (as in Figure~\ref{fig-fund-wedge}), consider a pin sequence starting on the two leftmost points of the alternation. At some stage, either the pin sequence attains length $2k$ (which we have assumed does not exist), or the pin has ``jumped'' a long way in the wedge alternation, giving the additional point required to form a wedge simple permutation of type 1 or 2 with length $2k$. Adjusting the bound $\lfloor n/2\rfloor < 2(16k^4)^{2k}$ to take into account that a wedge alternation of length $4k^2$ must be avoided gives the stated bound $n<2(256k^8)^{2k}$.
\end{proofsketch}

An immediate consequence of this decomposition is that within a suitably long simple permutation $\sigma$ we may find two ``almost disjoint'' simple subsequences by considering subsequences of the long proper pin sequence, parallel alternation or wedge simple permutation contained in $\sigma$. This result has consequences for certain permutation classes, which we will discuss in Subsection~\ref{subsec-finding-simples}.

\begin{corollary}[Brignall, Huczynska and Vatter~\cite{brignall:simple-permutat:a}]\label{sp2-main}
There is a function $f(k)$ such that every simple permutation of length at least $f(k)$ contains two simple subsequences, each of length at least $k$, sharing at most two entries.
\end{corollary}

To observe that it may be necessary that the two simple subsequences share exactly two entries, consider the family of type 2 wedge simple permutations (i.e.\ those of the form $m(2m)(m-1)(m+1)(m-2)(m+2)\cdots 1(2m-1)$), in which the first two entries are required in every simple subsequence.

\subsection{Simple Extensions}
In addition to finding the underlying structures in simple permutations, we may ask how, for an arbitrary permutation $\pi$, we may embed $\pi$ into a simple permutation, and how long this permutation has to be. The analogous problem for tournaments was solved by Erd\H os, Fried, Hajnal and Milner~\cite{erdhos:some-remarks-on:}, where they demonstrated that every tournament may be embedded in a simple tournament containing at most two extra vertices. Furthermore, in another paper published in the same year, Erd\H os, Hajnal and Milner~\cite{erdhos:simple-one-poin:} listed the cases when one point was not sufficient. For permutations, however, adding two points is rarely going to be sufficient. In fact, in the case of an increasing permutation $12\cdots n$ of length $n$, an additional $\lceil(n+1)/2\rceil$ points are both necessary and sufficient. The same bound can be obtained inductively for arbitrary permutations by means of the substitution decomposition, though it is not known when fewer points are sufficient.

\begin{theorem}[Brignall, Ru\v skuc and Vatter~\cite{brignall:embedding-combinat:}]
Every permutation $\pi$ on $n$ symbols has a simple extension with at most $\lceil(n+1)/2\rceil$ additional points.
\end{theorem}

Cases for graphs and posets are also obtainable by the same method. In the former, at most $\lceil\log_2(n+1)\rceil$ additional vertices are required, while for the latter, the bound is $\lceil(n+1)/2\rceil$, as in the permutation case.

%
%
%
%
%
%
%
%
%
%
\section{Permutation Classes with Finitely Many Simples}\label{sec-classes}
Considerable attention has been paid in recent years to the study of permutation classes which contain only finitely many simple permutations. Much is known about such classes both in terms of enumeration and structure, and of course these two considerations are not unrelated --- it is precisely the structure (given in terms of the substitution decomposition) of permutations lying in classes containing only finitely many simple permutations that gives these classes this well behaved enumeration.

\subsection{Substitution Closures}
As a preliminary to the results of this section, we consider permutation classes which may be described exactly by their (not necessarily finite) set of simple permutations. A class $\C$ of permutations is {\it substitution-closed\/} if $\sigma[\alpha_1,\dots,\alpha_m]\in\C$ for all $\sigma,\alpha_1,\dots,\alpha_m\in\C$.  The {\it substitution-closure\/}\footnote{The substitution closure has extensively been called the \emph{wreath closure} in the literature.} of a set $X$, $\wrc{X}$, is defined as the smallest substitution-closed class containing $X$.  (This concept is well-defined because the intersection of substitution-closed classes is substitution-closed, and the set of all permutations is substitution-closed.)

Letting $\Si(\C)$ denote the set of simple permutations in the class $\C$, we observe that $\Si(\C)=\Si(\wrc{\C})$: every permutation in $\C$ is an inflation of a member of $\Si(\C)$ so it follows (e.g., by induction) that $\C\subseteq\wrc{\Si(\C)}$.  Thus $\wrc{\C}\subseteq\wrc{\Si(C)}$, establishing that $\Si(\C)=\Si(\wrc{\C})$.  As substitution-closed classes are uniquely determined by their sets of simple permutations, $\wrc{\C}$ is the largest class with this property. For example, the substitution closure of $\Av(132)$ is the largest class whose only simple permutations are $1$, $12$, and $21$, which is precisely the class $\Av(2413,3142)$ of \emph{separable permutations}.

It is quite easy to decide if a permutation class given by a finite basis is substitution-closed:
\begin{proposition}[Albert and Atkinson~\cite{albert:simple-permutat:}]\label{decide-wreath-complete}
A permutation class is substitution-closed if and only if each of its basis elements is simple.
\end{proposition}

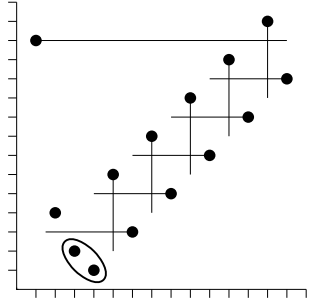
\begin{figure}
\begin{center}
\psset{xunit=0.01in, yunit=0.01in} \psset{linewidth=0.005in}
\begin{pspicture}(0,0)(155,155)
\psaxes[dy=10,dx=10,tickstyle=bottom,showorigin=false,labels=none](0,0)(150,150)
\psccurve[linewidth=0.01in](25,25)(30,10)(45,5)(40,20)
\psline(15,30)(60,30)
\psline(40,50)(80,50)
\psline(60,70)(100,70)
\psline(80,90)(120,90)
\psline(100,110)(140,110)
\psline(10,130)(140,130)
\psline(50,20)(50,60)
\psline(70,40)(70,80)
\psline(90,60)(90,100)
\psline(110,80)(110,120)
\psline(130,100)(130,140)
\pscircle*(10,130){0.03in}
\pscircle*(20,40){0.03in}
\pscircle*(30,20){0.03in}
\pscircle*(40,10){0.03in}
\pscircle*(50,60){0.03in}
\pscircle*(60,30){0.03in}
\pscircle*(70,80){0.03in}
\pscircle*(80,50){0.03in}
\pscircle*(90,100){0.03in}
\pscircle*(100,70){0.03in}
\pscircle*(110,120){0.03in}
\pscircle*(120,90){0.03in}
\pscircle*(130,140){0.03in}
\pscircle*(140,110){0.03in}
\end{pspicture}
\end{center}
\caption{A basis element of the substitution closure of $\Av(4321)$.}
\label{fig-inc-osc-antichain-variant}
\end{figure}

One may also wish to compute the basis of $\wrc{\C}$.  This is routine for classes with finitely many simple permutations (see Proposition~\ref{prop-wc-fin-simple}), but much less so in general.  An example of a finitely based class whose substitution-closure is infinitely based is $\Av(4321)$ --- the basis of its substitution-closure contains an infinite family of permutations, one of which is shown in Figure~\ref{fig-inc-osc-antichain-variant}.

The natural question is then:
\begin{question}\label{question-wreath-closure}
Given a finite basis $B$, is it decidable whether $\wrc{\Av(B)}$ is finitely based?\footnote{This question has recently been considered by Atkinson, Ru\v{s}kuc and Smith~\cite{atkinson:wreath-closed:}. The analogous question for graphs was raised by Giakoumakis~\cite{giakoumakis:on-the-closure-:} and has received a sizable amount of attention, see for example Zverovich~\cite{zverovich:a-finiteness-th:}.}
\end{question}

The importance of substitution-closed classes will become increasingly evident throughout this section. In particular, it is often fruitful to prove results about a class containing only finitely many simple permutations by first considering its substitution closure, and then adding basis restrictions to recover the original class.

\subsection{Algebraic Generating Functions}
When a class is enumerated by an algebraic generating function, we intuitively expect to find some recursive description of the permutations in the class. The converse intuition is also good, i.e.\ that classes constructed via some set of recursions should have algebraic generating functions. In a class $\C$ which has only finitely many simple permutations, any long permutation must be constructed recursively via the substitution decomposition, starting from the finite set $\Si(\C)$. Applying our intuition then immediately suggests the following result:

\begin{theorem}[Albert and Atkinson~\cite{albert:simple-permutat:}]\label{simple}
A permutation class with only finitely many simple permutations has a readily computable algebraic generating function.
\end{theorem}

Theorem~\ref{simple} was obtained by first proving it in the case of substitution-closed classes, and then demonstrating how adding an extra basis restriction does not affect the algebraicity of the generating function. A stronger result is true if, in addition to containing only finitely many simple permutations, our permutation class avoids a decreasing permutation of some length $n$, namely that such a class is enumerated by a rational generating function~\cite{albert:simple-permutat:}.

We may view permutation classes as subsets of the unique substitution-closed class containing the same set of simple permutations. When the substitution-closed class contains only finitely many simple permutations, Theorem~\ref{simple} shows that subclasses --- now thought of merely as particular subsets --- have algebraic generating functions. These subsets are not alone in this property, however; there are many other subsets which are enumerated by algebraic generating functions.

Calling any set of permutations $P$ a \emph{property}, we may say that a permutation $\pi$ \emph{satisfies} $P$ if $\pi\in P$. A set $\P$ of properties is said to be \emph{query-complete} if, for every simple permutation $\sigma$ and property $P\in\P$, one may determine whether $\sigma[\alpha_1,\ldots,\alpha_m]$ satisfies $P$ simply by knowing which properties of $\P$ each $\alpha_i$ satisfies. For example, avoidance of a given pattern is query-complete, since the property $\Av(\beta)$ lies in the query-complete set $\{\Av(\delta):\delta\leq\beta\}$. Moreover, since this query-complete set is finite, it is said to be \emph{finite query-complete}.

\begin{theorem}[Brignall, Huczynska and Vatter~\cite{brignall:simple-permutat:}]\label{sp1-really-main}
Let $\C$ be a permutation class containing only finitely many simple permutations, $\P$ a finite query-complete set of properties, and $\Q\subseteq\P$.  The generating function for the set of permutations in $\C$ satisfying every property in $\Q$ is algebraic over $\mathbb{Q}(x)$.
\end{theorem}

Properties known to lie in finite query-complete sets include the set of alternating permutations, even permutations, Dumont permutations of the first kind, and those avoiding any number of blocked or barred permutations. For example, a permutation $\pi$ of length $n$ is \emph{alternating} if, for all $i\in [2,n-1]$, $\pi(i)$ does not lie between $\pi(i-1)$ and $\pi(i+1)$. We can then explicitly state the finite query-complete set of properties:

\begin{lemma}
The set of properties consisting of
\begin{itemize}
\item $AL=\{\mbox{alternating permutations}\}$,
\item $BR=\{\mbox{permutations beginning with a rise, i.e., permutations with $\pi(1)<\pi(2)$}\}$,
\item $ER=\{\mbox{permutations ending with a rise}\}$, and
\item $\{1\}$.
\end{itemize}
is query-complete.
\end{lemma}
\begin{proof}
It is easy to show that $\{\{1\},BR,ER\}$ is query-complete:
\begin{eqnarray*}
\sigma[\alpha_1,\dots,\alpha_m]\in BR&\iff&
\alpha_1\in BR
\mbox{ or }
\left(\alpha_1=1\mbox{ and }\sigma\in BR\right),
\\
\sigma[\alpha_1,\dots,\alpha_m]\in ER&\iff&
\alpha_m\in ER
\mbox{ or }
\left(\alpha_m=1\mbox{ and }\sigma\in ER\right).
\end{eqnarray*}

If $\pi=\sigma[\alpha_1,\dots,\alpha_m]$ is an alternating permutation, we must have $\alpha_1,\dots,\alpha_m\in AL$. After checking that all the entries of $\pi$ up to and including the interval corresponding to $\sigma(i)$ are alternating, we must consider two cases: the first where $\sigma(i)>\sigma(i+1)$, the second $\sigma(i)<\sigma(i+1)$. We show only the former, the latter being analogous. Since $\sigma(i)>\sigma(i+1)$, $\pi$ contains a descent between its $\sigma(i)$ interval and its $\sigma(i+1)$ interval.  Thus $\alpha_i$ is allowed to be $1$ (i.e., $\alpha_i\in\{1\}$) only if $i=1$ or $\sigma(i-1)<\sigma(i)$, while if $\alpha_i\neq 1$ then we must have $\alpha_i\in ER$, and whether or not $\alpha_i$ is $1$ we must have $\alpha_{i+1}\in BR\cup\{1\}$.\end{proof}

The intersection of two properties known to lie in finite query-complete sets again lies in a finite query-complete set, and thus any combination of the properties listed above will lie in a finite query complete-set. It is also possible to show that the number of involutions of a permutation class containing only finitely many simple permutations is enumerated by an algebraic generating function, as is the cyclic closure of such a permutation class, and these may again be combined with other properties. Thus, for example, the number of alternating even involutions in a permutation class with only finitely many simple permutations is enumerated algebraically.

\subsection{Partial Well-Order}
We begin this subsection with a brief review of some of the known results about antichains and partial well-order. Most of these results apply to a variety of relational structures, though here we will restrict our attention to the permutation case.

An \emph{antichain} is a set of pairwise incomparable elements. Infinite antichains of permutations rely heavily on the structure of simple permutations to maintain their incomparability --- two permutations of an antichain are incomparable either because their skeletons are incomparable, or because there are incomparable blocks in the substitution decomposition. This latter situation suggests, heuristically, that we could instead have taken the smaller incomparable blocks as our antichain elements, and this motivates the definition of fundamental infinite antichains: an antichain $A$ is \emph{fundamental} if its \emph{closure}, defined by $\Sub(A)=\{\pi : \pi\leq\alpha \textrm{ for some }\alpha\in A\}$, contains no infinite antichains except for subsets of $A$ itself. The term ``fundamental'' is due to Murphy~\cite{murphy:restricted-perm:}; other authors (see, for example, Gustedt~\cite{gustedt:finiteness-theo:}) refer to such antichains as \emph{minimal}, because they are minimal under the following order on infinite antichains: $A \preceq B$ if $A$ is contained in the closure of $B$.

Any partial order is \emph{partially well-ordered} if it contains neither an infinite properly decreasing sequence nor an infinite antichain, though in the case of permutation classes we need only check the latter condition. Thus to show that a permutation class is not partially well ordered, it is sufficient to show that the class contains an infinite antichain. More strongly, it is sufficient to look for fundamental antichains --- the following result is obtainable using Nash-Williams' ``minimal bad sequence'' argument~\cite{nw:wqo}:

\begin{proposition}\label{fundamentals-exist} Every non-partially well-ordered permutation class contains a fundamental antichain.
\end{proposition}

On the other hand, proving that a permutation class is partially well-ordered is a somewhat harder task. The primary tool here is a result due to Higman~\cite{higman:ordering-by-div:}, which for brevity we will omit, as we do not require it.

Our intuition that simple permutations and fundamental antichains are connected is to some extent confirmed by studying classes containing only finitely many simple permutations --- we will now show that such classes are partially well-ordered. We first require the following straightforward observation:

\begin{proposition}\label{product-pwo}
The product $(P_1,\le_1)\times\cdots\times(P_s,\le_s)$ of a collection
of partial orders is partially well ordered if and only if each of them is partially well ordered.
\end{proposition}

No further results are needed for the required result. We will follow the proof given by Gustedt; Albert and Atkinson's proof uses Higman's Theorem, and is essentially the same in its approach in that it makes use of the substitution decomposition.

\begin{proposition}[Gustedt~\cite{gustedt:finiteness-theo:}; Albert and Atkinson~\cite{albert:simple-permutat:}]\label{fin-simples-pwo}
Every permutation class with only finitely many simple permutations is partially well-ordered.
\end{proposition}
\begin{proof}
Suppose to the contrary that the class $\C$ contains an infinite
antichain but only finitely many simple permutations.  By
Proposition~\ref{fundamentals-exist}, $\C$ contains an infinite fundamental
antichain. Moreover, there is an infinite subset $A$ of this antichain for which every element is
an inflation of the same simple permutation, say $\sigma$.  Let $\mathcal{D}$
denote the strict closure of $A$, noting that $A$ is
fundamental so $\mathcal{D}$ is partially well ordered.  It is easy to see that the permutation
containment order, when restricted to inflations of $\sigma$, is
isomorphic to a product order:
$\sigma[\alpha_1,\dots,\alpha_m]\le\sigma[\alpha_1',\dots,\alpha_m']$
if and only if $\alpha_i\le\alpha_i'$ for all $i\in[m]$.  However,
this implies that $A$ is an infinite antichain in a product
$\mathcal{D}\times\cdots\times\mathcal{D}$ of partially well-ordered posets, contradicting
Proposition~\ref{product-pwo}.
\end{proof}

\subsection{Finite Basis}
That a class containing only finitely many simple permutations is finitely based arises by first considering its substitution closure. Our first task is to compute the basis of a substitution closed class containing only finitely many simple permutations, which is easily done using Theorem~\ref{thm-schmerl-trotter}:

\begin{proposition}\label{prop-wc-fin-simple}
If the longest simple permutations in $\C$ have length $k$ then the basis elements of $\wrc{\C}$ have length at most $k+2$.
\end{proposition}
\begin{proof}
The basis of $\wrc{\C}$ is easily seen to consist of the minimal (under the pattern containment order) simple permutations not contained in $\C$ (cf.~Proposition~\ref{decide-wreath-complete}).  Let $\pi$ be such a permutation of length $n$.  Theorem~\ref{thm-schmerl-trotter} shows that $\pi$ contains a simple permutation $\sigma$ of length $n-1$ or $n-2$.  If $n\ge k+3$, then $\sigma\notin\C$, so $\sigma\notin\wrc{\C}$ and thus $\pi$ cannot lie in the basis of $\wrc{\C}$.
\end{proof}

For example, using this Proposition it can be computed that the wreath closure of $1$, $12$, $21$, and $2413$ is $\Av(3142,25314,246135,362514)$.

The finite basis result for arbitrary permutation classes containing only finitely many simple permutations now follows by recalling that all such classes are partially well-ordered. Thus:

\begin{theorem}[Murphy~\cite{murphy:restricted-perm:}; Albert and Atkinson~\cite{albert:simple-permutat:}]
Every permutation class containing only finitely many simple permutations is finitely based.
\end{theorem}

\begin{proof}
Let $\C$ be a class containing only finitely many permutations. By Proposition~\ref{prop-wc-fin-simple}, $\wrc{\C}$ is finitely based, and by Proposition~\ref{fin-simples-pwo} it is partially well ordered. The class $\C$ must therefore avoid all elements in the basis of $\wrc{\C}$, together with the minimal elements of $\wrc{\C}$ not belonging to $\C$, which form an antichain. By its partial well ordering any antichain in $\wrc{\C}$ is finite, and so there can only be finitely many basis elements of $\C$.
\end{proof}

\subsection{Finding Finitely Many Simples}\label{subsec-finding-simples}
In order to find whether the above properties apply to a given permutation class defined by its (finite) basis, it is of course necessary to know which simple permutations lie in the class, or at least whether there are only finitely many of them. In the first instance, we may use Theorem~\ref{thm-schmerl-trotter} simply to look for the simple permutations of each length $n=1,2,\ldots$, since if we encounter two consecutive integers for which the class has no simple permutations of those lengths, then there can be no simple permutations of any greater length.

One of the simplest classes with only finitely many simple permutations is $\Av(132)$, containing only $1$, $12$ and $21$. B\'ona~\cite{bona:the-number-of-p:} and Mansour and Vainshtein~\cite{mansour:counting-occurr:} showed that that for every $r$, the class of all permutations containing at most $r$ copies of $132$ has an algebraic generating function. The reason for this algebraicity is as we might hope: these classes contain only finitely many simple permutations, and this follows from a more general result that relies on Corollary~\ref{sp2-main}.

Denote by $\Av(\beta_1^{\le r_1},\beta_2^{\le r_2},\dots,\beta_k^{\le r_k})$
the set of permutations that have at most $r_1$ copies of $\beta_1$, $r_2$ copies of $\beta_2$, and so on. It should be clear that any such set forms a permutation class, although finding its basis is perhaps less obvious. Atkinson~\cite{atkinson:restricted-perm:} showed that the basis elements of the class can have length at most $\max\{(r_i+1)|\beta_i| : i\in[k]\}$. For example, $\Av(132^{\le 1})=\Av(1243$, $1342$, $1423$, $1432$, $2143$, $35142$, $354162$, $461325$, $465132)$. We then have the following result:

\begin{theorem}[Brignall, Huczynska and Vatter~\cite{brignall:simple-permutat:a}]\label{sp2-main-cor}
If the class $\Av(\beta_1,\beta_2,\dots,\beta_k)$ contains only finitely many simple permutations then the class $\Av(\beta_1^{\le r_1}$, $\beta_2^{\le r_2}$, $\dots$, $\beta_k^{\le r_k})$ also contains only finitely many simple permutations for all choices of nonnegative integers $r_1$, $r_2$, $\dots$, $r_k$.
\end{theorem}
\begin{proof}
We need to show that for any choice of nonnegative integers $r_1,r_2,\dots,r_k$, only finitely many simple permutations contain at most $r_i$ copies of $\beta_i$ for each $i\in[k]$.  We may suppose that $|\beta_i|\ge 3$ for all $i\in [k]$ since if any $\beta_i$ is of length $1$ or $2$ then the theorem follows easily.  We now proceed by induction. The base case, arising when $r_i=0$ for all $i$, follows trivially, so suppose that some $r_j>0$ and set
\[
g(r_1,r_2,\dots,r_k)=f(g(r_1,r_2,\dots,r_{j-1},\lfloor r_j/2\rfloor, r_{j+1},\dots,r_k)),
\]
where $f$ is the function from Corollary~\ref{sp2-main}.  By that result, every simple permutation $\pi$ of length at least $g(r_1,r_2,\dots,r_k)$ contains two simple subsequences of length at least $f(g(r_1,r_2,\dots,r_{j-1},\lfloor r_j/2\rfloor, r_{j+1},\dots,r_k))$, and by induction each of these simple subsequences contains more than $\lfloor r_j/2 \rfloor $ copies of $\beta_j$.  Moreover, because these simple subsequences share at most two entries, their copies of $\beta_j$ are distinct, and thus $\pi$ contains more than $r_j$ copies of $\beta_j$, proving the theorem.
\end{proof}

For a general answer to the decidability question, we turn to Theorem~\ref{sp2-really-main}, which reduces the task to checking whether the permutation class in question contains arbitrarily long proper pin sequences, parallel alternations or wedge simple alternations of types $1$ and $2$. This may be done algorithmically:

\begin{theorem}[Brignall, Ru\v{s}kuc and Vatter~\cite{brignall:simple-permutat:b}]\label{sp3-main}
It is possible to decide if a permutation class given by a finite basis contains infinitely many simple permutations.
\end{theorem}

\begin{proofsketch} To determine whether a permutation class $\Av(B)$ contains only finitely many parallel alternations oriented $\backslash\backslash$, we need only check that there is a permutation in $B$ that is contained in such an alternation. More simply, it is sufficient to verify that $B$ contains a permutation in the class $\Av(123,2413,3412)$. To check the other orientations, we need simply test the same condition for every symmetry of this class. Similarly, to check the wedge simple permutations of types $1$ and $2$, one needs to ensure that $B$ contains an element of every symmetry of $\Av(1243$, $1324$, $1423$, $1432$, $2431$, $3124$, $4123$, $4132$, $4231$, $4312)$ and $\Av(2134$, $ 2143$, $3124$, $3142$, $3241$, $3412$, $4123$, $4132$, $4231$, $4312)$ respectively.

Thus it remains to determine whether $\Av(B)$ contains arbitrarily long proper pin sequences. This is done by encoding proper pin sequences as ``strict pin words'' over the four-letter alphabet of \emph{directions} $\{L,R,U,D\}$ and subsequences of proper pin sequences as pin words over an eight-letter alphabet consisting of the four directions and four \emph{numerals} $\{1,2,3,4\}$. These numerals correspond to the four quadrants around an origin which is placed ``close'' to the first two points of the pin sequence in such a way that whenever a pin is not included in the subsequence the next pin that is included is encoded as a numeral corresponding to the quadrant in which it lies. The permutation containment order restricted to pin sequences and their subsequences corresponds to an order on these words.

For every $\beta\in B$, list all (non-strict) pin words corresponding to $\beta$. It may be shown that the set of strict pin words containing any non-strict pin word forms a regular language, and hence the union of all strict pin words containing any pin word corresponding to some $\beta$ forms a regular language. The complement of this set of strict pin words corresponds to the proper pin sequences lying in $\Av(B)$, and also forms a regular language. It is decidable whether a regular language contains arbitrarily long words or not, from which the result follows.
\end{proofsketch}

\subsection{Algorithms}
\paragraph{Linear Time Membership.}
Bose, Buss and Lubiw~\cite{bose:matching} showed that deciding whether a given permutation lies in some permutation class is in general NP-complete, but that one may use the substitution decomposition to decide whether a permutation is separable in polynomial time. Out of some of the recent machinery surveyed here comes an indication that, given a permutation class $\C$ that contains only finitely many simple permutations, it may be decided in linear time whether an arbitrary permutation $\pi$ of length $n$ lies in $\C$. The approach relies first and foremost on the fact that we may compute the substitution decomposition of any permutation in linear time, as mentioned in Subsection~\ref{subsec-subst-decomp}. We begin by first performing some precomputations specific to the class $\C$, all of which may be done essentially in constant time:
\begin{itemize}
\item Compute $\Si(\C)$, the set of simple permutations in $\C$.
\item Compute the basis $B$ of $\C$, noting that permutations in $B$ can be no longer than $\displaystyle \max_{\sigma\in\Si(\C)}|\sigma| + 2$ by the Schmerl-Trotter Theorem~\ref{thm-schmerl-trotter}.
\item For every $\beta$ either lying in $B$ or contained in a permutation lying in $B$, list all expressions of $\beta$ as a \emph{lenient inflation} (an inflation $\sigma[\gamma_1,\ldots,\gamma_m]$ in which the $\gamma_i$'s are allowed to be empty) of each $\sigma\in\Si(\C)$.
\end{itemize}

With this one-time work done, we now take our candidate permutation $\pi$ of length $n$ and compute its substitution decomposition, $\pi=\sigma[\alpha_1,\ldots,\alpha_m]$. Now, after first trivially checking that the skeleton $\sigma$ lies in $\C$, we look at all the expressions of each $\beta\in B$ as lenient inflations of $\sigma$. Note that if some $\beta\leq\pi$, there must exist an expression of $\beta$ as a lenient inflation $\beta=\sigma[\gamma_1,\ldots,\gamma_m]$ so that $\gamma_i\leq\alpha_i$ for every $i=1,\ldots,m$.

Thus, taking each lenient inflation $\beta=\sigma[\gamma_1,\ldots,\gamma_m]$ in turn, we look recursively at each block, testing to see if $\gamma_i\leq\alpha_i$ is true. Though this recursion makes the linear-time complexity non-obvious, note that the number of levels of recursion that are required cannot be more than the maximum depth of the substitution decomposition tree, which itself cannot have more than $2n$ nodes. The recursion will eventually reduce the problem to making only trivial comparisons, each of which is immediately answerable in constant time.

\paragraph{Longest Common Pattern.}
A possible generalisation of the pattern containment problem is that of computing the longest common pattern of two permutations. Bouvel and Rossin~\cite{bouvel:longest-common-:} demonstrate a general algorithm to find the longest common pattern between two permutations $\pi_1$ and $\pi_2$, relying on the substitution decomposition of either $\pi_1$ or $\pi_2$. For general $\pi_1$ and $\pi_2$ this does not run in polynomial time, but in the special case where $\pi_1$ comes from a permutation class whose simple permutations are of length at most $d$, their algorithm runs in $O(\min(n_1,n_2)n_1n_2^{2d+2})$, where $n_1=|\pi_1|$ and $n_2=|\pi_2|$. In particular, they give an explicit $O(n^8)$ algorithm for finding the longest common pattern between two permutations of length $n$ given one is separable.

\section{Concluding Remarks}

As we have seen, much is known about permutation classes containing only finitely many simple permutations. On the other hand, little is known in general about classes containing infinitely many simple permutations. The way in which many of the results presented in Section~\ref{sec-classes} are obtained suggests that a first step would be to restrict our attention to substitution-closed classes. There is one result to support this approach in the more general context: it is known that a substitution-closed class has an algebraic generating function if and only if the simple permutations of the class are also enumerated by an algebraic generating function~\cite{albert:simple-permutat:}.

Related to this is the question of partial well-order --- we have seen that simple permutations are intimately related to the structure of fundamental antichains, but precisely how they are related remains unknown. In particular, if a permutation class is partially well-ordered, what can be said about its set of simple permutations?

Recalling that the substitution decomposition is defined analogously for all relational structures, there is no reason why some of the results reviewed here cannot be extended to other relational structures. While questions such as the enumeration of hereditary properties are not obviously extendable (since isomorphism between structures must be taken into account), questions relating to the decomposition of simple or indecomposable objects and to the study of partial well-order in hereditary properties are likely to behave in a similar way for all structures.

\paragraph{Acknowledgments.} The author thanks Vince Vatter and the anonymous referee for their helpful and enriching comments. Thanks is also owed to Mathilde Bouvel for fruitful discussions.

\bibliographystyle{abbrv}
\bibliography{../refs}

\end{document}